\documentclass{amsart}
\usepackage{amscd}
\let\cal\mathcal
\def\Ascr{{\cal A}}
\def\Bscr{{\cal B}}
\def\Cscr{{\cal C}}
\def\Dscr{{\cal D}}

\def\Fscr{{\cal F}}

\def\Oscr{{\cal O}}

\def\Sscr{{\cal S}}
\def\Tscr{{\cal T}}

\let\blb\mathbb

\def \ZZ{{\blb Z}}

\def\Id{\operatorname{id}}

\def\Lotimes{\overset{L}{\otimes}}

\def\coh{\mathop{\text{\upshape{coh}}}}

\def\Ext{\operatorname {Ext}}
\def\Hom{\operatorname {Hom}}
\def\RHom{\operatorname {RHom}}

\def\cd{\operatorname {cd}}

\def\coker{\operatorname {coker}}
\def\ker{\operatorname {ker}}

\def\r{\rightarrow}

\newtheorem{lemma}{Lemma}[section]

\newtheorem{theorem}[lemma]{Theorem}

\theoremstyle{definition}

\newtheorem{question}[lemma]{Question}

{

}

\theoremstyle{remark}

\newtheorem{remark}[lemma]{Remark}

\newdimen\uboxsep \uboxsep=1ex
\def\uboxn#1{\vtop to 0pt{\hrule height 0pt depth 0pt\vskip\uboxsep
\hbox to 0pt{\hss #1\hss}\vss}}

\def\uboxs#1{\vbox to 0pt{\vss\hbox to 0pt{\hss #1\hss}
\vskip\uboxsep\hrule height 0pt depth 0pt}}

\makeatletter

\def\diagram{\m@th\leftwidth=\z@ \rightwidth=\z@ \topheight=\z@
\botheight=\z@ \setbox\@picbox\hbox\bgroup}

\def\enddiagram{\egroup\wd\@picbox\rightwidth\unitlength
\ht\@picbox\topheight\unitlength \dp\@picbox\botheight\unitlength
\hskip\leftwidth\unitlength\box\@picbox}

\def\bfig{\begin{diagram}}
\def\efig{\end{diagram}}
\newcount\wideness \newcount\leftwidth \newcount\rightwidth
\newcount\highness \newcount\topheight \newcount\botheight

\def\ratchet#1#2{\ifnum#1<#2 \global #1=#2 \fi}

\def\putbox(#1,#2)#3{%
\horsize{\wideness}{#3} \divide\wideness by 2
{\advance\wideness by #1 \ratchet{\rightwidth}{\wideness}}
{\advance\wideness by -#1 \ratchet{\leftwidth}{\wideness}}
\vertsize{\highness}{#3} \divide\highness by 2
{\advance\highness by #2 \ratchet{\topheight}{\highness}}
{\advance\highness by -#2 \ratchet{\botheight}{\highness}}
\put(#1,#2){\makebox(0,0){$#3$}}}

\def\putlbox(#1,#2)#3{%
\horsize{\wideness}{#3}
{\advance\wideness by #1 \ratchet{\rightwidth}{\wideness}}
{\ratchet{\leftwidth}{-#1}}
\vertsize{\highness}{#3} \divide\highness by 2
{\advance\highness by #2 \ratchet{\topheight}{\highness}}
{\advance\highness by -#2 \ratchet{\botheight}{\highness}}
\put(#1,#2){\makebox(0,0)[l]{$#3$}}}

\def\putrbox(#1,#2)#3{%
\horsize{\wideness}{#3}
{\ratchet{\rightwidth}{#1}}
{\advance\wideness by -#1 \ratchet{\leftwidth}{\wideness}}
\vertsize{\highness}{#3} \divide\highness by 2
{\advance\highness by #2 \ratchet{\topheight}{\highness}}
{\advance\highness by -#2 \ratchet{\botheight}{\highness}}
\put(#1,#2){\makebox(0,0)[r]{$#3$}}}

\def\adjust[#1]{} 

\newcount \coefa
\newcount \coefb
\newcount \coefc
\newcount\tempcounta
\newcount\tempcountb
\newcount\tempcountc
\newcount\tempcountd
\newcount\xext
\newcount\yext
\newcount\xoff
\newcount\yoff
\newcount\gap%
\newcount\arrowtypea
\newcount\arrowtypeb
\newcount\arrowtypec
\newcount\arrowtyped
\newcount\arrowtypee
\newcount\height
\newcount\width
\newcount\xpos
\newcount\ypos
\newcount\run
\newcount\rise
\newcount\arrowlength
\newcount\halflength
\newcount\arrowtype
\newdimen\tempdimen
\newdimen\xlen
\newdimen\ylen
\newsavebox{\tempboxa}%
\newsavebox{\tempboxb}%
\newsavebox{\tempboxc}%

\newdimen\w@dth

\def\setw@dth#1#2{\setbox\z@\hbox{$#1$}\w@dth=\wd\z@
\setbox\@ne\hbox{$#2$}\ifnum\w@dth<\wd\@ne \w@dth=\wd\@ne \fi
\advance\w@dth by 1.2em}

\def\t@^#1_#2{\allowbreak\def\n@one{#1}\def\n@two{#2}\mathrel
{\setw@dth{#1}{#2}
\mathop{\hbox to \w@dth{\rightarrowfill}}\limits
\ifx\n@one\empty\else ^{\box\z@}\fi
\ifx\n@two\empty\else _{\box\@ne}\fi}}
\def\t@@^#1{\@ifnextchar_ {\t@^{#1}}{\t@^{#1}_{}}}
\def\to{\@ifnextchar^ {\t@@}{\t@@^{}}}

\def\t@left^#1_#2{\def\n@one{#1}\def\n@two{#2}\mathrel{\setw@dth{#1}{#2}
\mathop{\hbox to \w@dth{\leftarrowfill}}\limits
\ifx\n@one\empty\else ^{\box\z@}\fi
\ifx\n@two\empty\else _{\box\@ne}\fi}}
\def\t@@left^#1{\@ifnextchar_ {\t@left^{#1}}{\t@left^{#1}_{}}}
\def\toleft{\@ifnextchar^ {\t@@left}{\t@@left^{}}}

\def\two@^#1_#2{\def\n@one{#1}\def\n@two{#2}\mathrel{\setw@dth{#1}{#2}
\mathop{\vcenter{\hbox to \w@dth{\rightarrowfill}\kern-1.7ex
                 \hbox to \w@dth{\rightarrowfill}}%
       }\limits
\ifx\n@one\empty\else ^{\box\z@}\fi
\ifx\n@two\empty\else _{\box\@ne}\fi}}
\def\tw@@^#1{\@ifnextchar_ {\two@^{#1}}{\two@^{#1}_{}}}
\def\two{\@ifnextchar^ {\tw@@}{\tw@@^{}}}

\def\tofr@^#1_#2{\def\n@one{#1}\def\n@two{#2}\mathrel{\setw@dth{#1}{#2}
\mathop{\vcenter{\hbox to \w@dth{\rightarrowfill}\kern-1.7ex
                 \hbox to \w@dth{\leftarrowfill}}%
       }\limits
\ifx\n@one\empty\else ^{\box\z@}\fi
\ifx\n@two\empty\else _{\box\@ne}\fi}}
\def\t@fr@^#1{\@ifnextchar_ {\tofr@^{#1}}{\tofr@^{#1}_{}}}
\def\tofro{\@ifnextchar^ {\t@fr@}{\t@fr@^{}}}

\def\mon{\mathop{\m@th\hbox to
      14.6\P@{\lasyb\char'51\hskip-2.1\P@$\arrext$\hss
$\mathord\rightarrow$}}\limits} 
\def\leftmono{\mathrel{\m@th\hbox to
14.6\P@{$\mathord\leftarrow$\hss$\arrext$\hskip-2.1\P@\lasyb\char'50%
}}\limits} 
\mathchardef\arrext="0200       

\def\settypes(#1,#2,#3){\arrowtypea#1 \arrowtypeb#2 \arrowtypec#3}
\def\settoheight#1#2{\setbox\@tempboxa\hbox{#2}#1\ht\@tempboxa\relax}%
\def\settodepth#1#2{\setbox\@tempboxa\hbox{#2}#1\dp\@tempboxa\relax}%
\def\settokens[#1`#2`#3`#4]{%
     \def\tokena{#1}\def\tokenb{#2}\def\tokenc{#3}\def\tokend{#4}}
\def\setsqparms[#1`#2`#3`#4;#5`#6]{%
\arrowtypea #1
\arrowtypeb #2
\arrowtypec #3
\arrowtyped #4
\width #5
\height #6
}
\def\setpos(#1,#2){\xpos=#1 \ypos#2}

\def\settriparms[#1`#2`#3;#4]{\settripairparms[#1`#2`#3`1`1;#4]}%

\def\settripairparms[#1`#2`#3`#4`#5;#6]{%
\arrowtypea #1
\arrowtypeb #2
\arrowtypec #3
\arrowtyped #4
\arrowtypee #5
\width #6
\height #6
}

\def\resetparms{\settripairparms[1`1`1`1`1;500]\width 500}

\resetparms

\def\mvector(#1,#2)#3{
\put(0,0){\vector(#1,#2){#3}}%
\put(0,0){\vector(#1,#2){26}}%
}
\def\evector(#1,#2)#3{{
\arrowlength #3
\put(0,0){\vector(#1,#2){\arrowlength}}%
\advance \arrowlength by-30
\put(0,0){\vector(#1,#2){\arrowlength}}%
}}

\def\horsize#1#2{%
\settowidth{\tempdimen}{$#2$}%
#1=\tempdimen
\divide #1 by\unitlength
}

\def\vertsize#1#2{%
\settoheight{\tempdimen}{$#2$}%
#1=\tempdimen
\settodepth{\tempdimen}{$#2$}%
\advance #1 by\tempdimen
\divide #1 by\unitlength
}

\def\putvector(#1,#2)(#3,#4)#5#6{{%
\ifnum3<\arrowtype
\putdashvector(#1,#2)(#3,#4)#5\arrowtype
\else
\ifnum\arrowtype<-3
\putdashvector(#1,#2)(#3,#4)#5\arrowtype
\else
\xpos=#1
\ypos=#2
\run=#3
\rise=#4
\arrowlength=#5
\ifnum \arrowtype<0
    \ifnum \run=0
        \advance \ypos by-\arrowlength
    \else
        \tempcounta \arrowlength
        \multiply \tempcounta by\rise
        \divide \tempcounta by\run
        \ifnum\run>0
            \advance \xpos by\arrowlength
            \advance \ypos by\tempcounta
        \else
            \advance \xpos by-\arrowlength
            \advance \ypos by-\tempcounta
        \fi
    \fi
    \multiply \arrowtype by-1
    \multiply \rise by-1
    \multiply \run by-1
\fi
\ifcase \arrowtype
\or \put(\xpos,\ypos){\vector(\run,\rise){\arrowlength}}%
\or \put(\xpos,\ypos){\mvector(\run,\rise)\arrowlength}%
\or \put(\xpos,\ypos){\evector(\run,\rise){\arrowlength}}%
\fi\fi\fi
}}

\def\putsplitvector(#1,#2)#3#4{
\xpos #1
\ypos #2
\arrowtype #4
\halflength #3
\arrowlength #3
\gap 140
\advance \halflength by-\gap
\divide \halflength by2
\ifnum\arrowtype>0
   \ifcase \arrowtype
   \or \put(\xpos,\ypos){\line(0,-1){\halflength}}%
       \advance\ypos by-\halflength
       \advance\ypos by-\gap
       \put(\xpos,\ypos){\vector(0,-1){\halflength}}%
   \or \put(\xpos,\ypos){\line(0,-1)\halflength}%
       \put(\xpos,\ypos){\vector(0,-1)3}%
       \advance\ypos by-\halflength
       \advance\ypos by-\gap
       \put(\xpos,\ypos){\vector(0,-1){\halflength}}%
   \or \put(\xpos,\ypos){\line(0,-1)\halflength}%
       \advance\ypos by-\halflength
       \advance\ypos by-\gap
       \put(\xpos,\ypos){\evector(0,-1){\halflength}}%
   \fi
\else \arrowtype=-\arrowtype
   \ifcase\arrowtype
   \or \advance \ypos by-\arrowlength
       \put(\xpos,\ypos){\line(0,1){\halflength}}%
       \advance\ypos by\halflength
       \advance\ypos by\gap
       \put(\xpos,\ypos){\vector(0,1){\halflength}}%
   \or \advance \ypos by-\arrowlength
       \put(\xpos,\ypos){\line(0,1)\halflength}%
       \put(\xpos,\ypos){\vector(0,1)3}%
       \advance\ypos by\halflength
       \advance\ypos by\gap
       \put(\xpos,\ypos){\vector(0,1){\halflength}}%
   \or \advance \ypos by-\arrowlength
       \put(\xpos,\ypos){\line(0,1)\halflength}%
       \advance\ypos by\halflength
       \advance\ypos by\gap
       \put(\xpos,\ypos){\evector(0,1){\halflength}}%
   \fi
\fi
}

\def\putmorphism(#1)(#2,#3)[#4`#5`#6]#7#8#9{{%
\run #2
\rise #3
\ifnum\rise=0
  \puthmorphism(#1)[#4`#5`#6]{#7}{#8}#9%
\else\ifnum\run=0
  \putvmorphism(#1)[#4`#5`#6]{#7}{#8}#9%
\else
\setpos(#1)%
\arrowlength #7
\arrowtype #8
\ifnum\run=0
\else\ifnum\rise=0
\else
\ifnum\run>0
    \coefa=1
\else
   \coefa=-1
\fi
\ifnum\arrowtype>0
   \coefb=0
   \coefc=-1
\else
   \coefb=\coefa
   \coefc=1
   \arrowtype=-\arrowtype
\fi
\width=2
\multiply \width by\run
\divide \width by\rise
\ifnum \width<0  \width=-\width\fi
\advance\width by60
\if l#9 \width=-\width\fi
\putbox(\xpos,\ypos){#4}
{\multiply \coefa by\arrowlength
\advance\xpos by\coefa
\multiply \coefa by\rise
\divide \coefa by\run
\advance \ypos by\coefa
\putbox(\xpos,\ypos){#5} }%
{\multiply \coefa by\arrowlength
\divide \coefa by2
\advance \xpos by\coefa
\advance \xpos by\width
\multiply \coefa by\rise
\divide \coefa by\run
\advance \ypos by\coefa
\if l#9%
   \putrbox(\xpos,\ypos){#6}%
\else\if r#9%
   \putlbox(\xpos,\ypos){#6}%
\fi\fi }%
{\multiply \rise by-\coefc
\multiply \run by-\coefc
\multiply \coefb by\arrowlength
\advance \xpos by\coefb
\multiply \coefb by\rise
\divide \coefb by\run
\advance \ypos by\coefb
\multiply \coefc by70
\advance \ypos by\coefc
\multiply \coefc by\run
\divide \coefc by\rise
\advance \xpos by\coefc
\multiply \coefa by140
\multiply \coefa by\run
\divide \coefa by\rise
\advance \arrowlength by\coefa
\ifcase\arrowtype
\or \put(\xpos,\ypos){\vector(\run,\rise){\arrowlength}}%
\or \put(\xpos,\ypos){\mvector(\run,\rise){\arrowlength}}%
\or \put(\xpos,\ypos){\evector(\run,\rise){\arrowlength}}%
\fi}\fi\fi\fi\fi}}

\newcount\numbdashes \newcount\lengthdash \newcount\increment

\def\howmanydashes{
\numbdashes=\arrowlength \lengthdash=40
\divide\numbdashes by \lengthdash
\lengthdash=\arrowlength
\divide\lengthdash by \numbdashes
\increment=\lengthdash
\multiply\lengthdash by 3
\divide\lengthdash by 5
}

\def\putdashvector(#1)(#2,#3)#4#5{%
\ifnum#3=0 \putdashhvector(#1){#4}#5
\else
\ifnum#2=0
\putdashvvector(#1){#4}#5\fi\fi}

\def\putdashhvector(#1,#2)#3#4{{%
\arrowlength=#3 \howmanydashes
\multiput(#1,#2)(\increment,0){\numbdashes}%
{\vrule height .4pt width \lengthdash\unitlength}
\arrowtype=#4 \xpos=#1
\ifnum\arrowtype<0 \advance\arrowtype by 7 \fi
\ifcase\arrowtype
\or \advance\xpos by 10
    \put(\xpos,#2){\vector(-1,0){\lengthdash}}
    \advance\xpos by 40
    \put(\xpos,#2){\vector(-1,0){\lengthdash}}
\or \advance \xpos by 10
    \put(\xpos,#2){\vector(-1,0){\lengthdash}}
    \advance\xpos by  \arrowlength
    \advance\xpos by  -50
    \put(\xpos,#2){\vector(-1,0){\lengthdash}}
\or \advance\xpos by 10
    \put(\xpos,#2){\vector(-1,0){\lengthdash}}
\or \advance\xpos by \arrowlength
    \advance\xpos by -\lengthdash
    \put(\xpos,#2){\vector(1,0){\lengthdash}}
\or {\advance\xpos by 10
    \put(\xpos,#2){\vector(1,0){\lengthdash}}}
    \advance\xpos by \arrowlength
    \advance\xpos by -\lengthdash
    \put(\xpos,#2){\vector(1,0){\lengthdash}}
\or \advance\xpos by \arrowlength
    \advance\xpos by -\lengthdash
    \put(\xpos,#2){\vector(1,0){\lengthdash}}
    \advance\xpos by -40
    \put(\xpos,#2){\vector(1,0){\lengthdash}}
   \fi
}}

\def\putdashvvector(#1,#2)#3#4{{%
\arrowlength=#3 \howmanydashes
\ypos=#2 \advance\ypos by -\arrowlength
\multiput(#1,#2)(0,\increment){\numbdashes}%
    {\vrule width .4pt height \lengthdash\unitlength}
\arrowtype=#4 \ypos=#2
\ifnum\arrowtype<0 \advance\arrowtype by 7 \fi
\ifcase\arrowtype
\or \advance\ypos by \arrowlength \advance\ypos by -40
    \put(#1,\ypos){\vector(0,1){\lengthdash}}
    \advance\ypos by -40
    \put(#1,\ypos){\vector(0,1){\lengthdash}}
\or \advance\ypos by 10
    \put(#1,\ypos){\vector(0,1){\lengthdash}}
    \advance\ypos by \arrowlength \advance\ypos by -40
    \put(#1,\ypos){\vector(0,1){\lengthdash}}
\or \advance\ypos by \arrowlength \advance\ypos by -40
    \put(#1,\ypos){\vector(0,1){\lengthdash}}
\or \advance\ypos by 10
    \put(#1,\ypos){\vector(0,-1){\lengthdash}}
\or \advance\ypos by 10
    \put(#1,\ypos){\vector(0,-1){\lengthdash}}
    \advance\ypos by \arrowlength \advance\ypos by -40
    \put(#1,\ypos){\vector(0,-1){\lengthdash}}
\or \advance\ypos by 10
    \put(#1,\ypos){\vector(0,-1){\lengthdash}}
    \advance\ypos by 40
    \put(#1,\ypos){\vector(0,-1){\lengthdash}}
\fi
}}

\def\puthmorphism(#1,#2)[#3`#4`#5]#6#7#8{{%
\xpos #1
\ypos #2
\width #6
\arrowlength #6
\arrowtype=#7
\putbox(\xpos,\ypos){#3\vphantom{#4}}%
{\advance \xpos by\arrowlength
\putbox(\xpos,\ypos){\vphantom{#3}#4}}%
\horsize{\tempcounta}{#3}%
\horsize{\tempcountb}{#4}%
\divide \tempcounta by2
\divide \tempcountb by2
\advance \tempcounta by30
\advance \tempcountb by30
\advance \xpos by\tempcounta
\advance \arrowlength by-\tempcounta
\advance \arrowlength by-\tempcountb
\putvector(\xpos,\ypos)(1,0)\arrowlength\arrowtype
\divide \arrowlength by2
\advance \xpos by\arrowlength
\vertsize{\tempcounta}{#5}%
\divide\tempcounta by2
\advance \tempcounta by20
\if a#8 %
   \advance \ypos by\tempcounta
   \putbox(\xpos,\ypos){#5}%
\else
   \advance \ypos by-\tempcounta
   \putbox(\xpos,\ypos){#5}%
\fi}}

\def\putvmorphism(#1,#2)[#3`#4`#5]#6#7#8{{%
\xpos #1
\ypos #2
\arrowlength #6
\arrowtype #7
\settowidth{\xlen}{$#5$}%
\putbox(\xpos,\ypos){#3}%
{\advance \ypos by-\arrowlength
\putbox(\xpos,\ypos){#4}}%
{\advance\arrowlength by-140
\advance \ypos by-70
\ifdim\xlen>0pt
   \if m#8%
      \putsplitvector(\xpos,\ypos)\arrowlength\arrowtype
   \else
   \putvector(\xpos,\ypos)(0,-1)\arrowlength\arrowtype
   \fi
\else
   \putvector(\xpos,\ypos)(0,-1)\arrowlength\arrowtype
\fi}%
\ifdim\xlen>0pt
   \divide \arrowlength by2
   \advance\ypos by-\arrowlength
   \if l#8%
      \advance \xpos by-40
      \putrbox(\xpos,\ypos){#5}%
   \else\if r#8%
      \advance \xpos by40
      \putlbox(\xpos,\ypos){#5}%
   \else
      \putbox(\xpos,\ypos){#5}%
   \fi\fi
\fi
}}

\def\putsquarep<#1>(#2)[#3;#4`#5`#6`#7]{{%
\setsqparms[#1]%
\setpos(#2)%
\settokens[#3]%
\puthmorphism(\xpos,\ypos)[\tokenc`\tokend`{#7}]{\width}{\arrowtyped}b%
\advance\ypos by \height
\puthmorphism(\xpos,\ypos)[\tokena`\tokenb`{#4}]{\width}{\arrowtypea}a%
\putvmorphism(\xpos,\ypos)[``{#5}]{\height}{\arrowtypeb}l%
\advance\xpos by \width
\putvmorphism(\xpos,\ypos)[``{#6}]{\height}{\arrowtypec}r%
}}

\def\putsquare{\@ifnextchar <{\putsquarep}{\putsquarep%
   <\arrowtypea`\arrowtypeb`\arrowtypec`\arrowtyped;\width`\height>}}
\def\square{\@ifnextchar< {\squarep}{\squarep
   <\arrowtypea`\arrowtypeb`\arrowtypec`\arrowtyped;\width`\height>}}
\def\squarep<#1>[#2`#3`#4`#5;#6`#7`#8`#9]{{
\setsqparms[#1]
\diagram
\putsquarep<\arrowtypea`\arrowtypeb`\arrowtypec`
\arrowtyped;\width`\height>
(0,0)[#2`#3`#4`{#5};#6`#7`#8`{#9}]
\enddiagram
}}                                                 
\def\putptrianglep<#1>(#2,#3)[#4`#5`#6;#7`#8`#9]{{%
\settriparms[#1]%
\xpos=#2 \ypos=#3
\advance\ypos by \height
\puthmorphism(\xpos,\ypos)[#4`#5`{#7}]{\height}{\arrowtypea}a%
\putvmorphism(\xpos,\ypos)[`#6`{#8}]{\height}{\arrowtypeb}l%
\advance\xpos by\height
\putmorphism(\xpos,\ypos)(-1,-1)[``{#9}]{\height}{\arrowtypec}r%
}}

\def\putptriangle{\@ifnextchar <{\putptrianglep}{\putptrianglep
   <\arrowtypea`\arrowtypeb`\arrowtypec;\height>}}
\def\ptriangle{\@ifnextchar <{\ptrianglep}{\ptrianglep
   <\arrowtypea`\arrowtypeb`\arrowtypec;\height>}}
\def\ptrianglep<#1>[#2`#3`#4;#5`#6`#7]{{
\settriparms[#1]
\diagram
\putptrianglep<\arrowtypea`\arrowtypeb`
\arrowtypec;\height>
(0,0)[#2`#3`#4;#5`#6`{#7}]
\enddiagram
}}                                            

\def\putqtrianglep<#1>(#2,#3)[#4`#5`#6;#7`#8`#9]{{%
\settriparms[#1]%
\xpos=#2 \ypos=#3
\advance\ypos by\height
\puthmorphism(\xpos,\ypos)[#4`#5`{#7}]{\height}{\arrowtypea}a%
\putmorphism(\xpos,\ypos)(1,-1)[``{#8}]{\height}{\arrowtypeb}l%
\advance\xpos by\height
\putvmorphism(\xpos,\ypos)[`#6`{#9}]{\height}{\arrowtypec}r%
}}

\def\putqtriangle{\@ifnextchar <{\putqtrianglep}{\putqtrianglep
   <\arrowtypea`\arrowtypeb`\arrowtypec;\height>}}
\def\qtriangle{\@ifnextchar <{\qtrianglep}{\qtrianglep
   <\arrowtypea`\arrowtypeb`\arrowtypec;\height>}}
\def\qtrianglep<#1>[#2`#3`#4;#5`#6`#7]{{
\settriparms[#1]
\width=\height                                
\diagram
\putqtrianglep<\arrowtypea`\arrowtypeb`
\arrowtypec;\height>
(0,0)[#2`#3`#4;#5`#6`{#7}]
\enddiagram
}}

\def\putdtrianglep<#1>(#2,#3)[#4`#5`#6;#7`#8`#9]{{%
\settriparms[#1]%
\xpos=#2 \ypos=#3
\puthmorphism(\xpos,\ypos)[#5`#6`{#9}]{\height}{\arrowtypec}b%
\advance\xpos by \height \advance\ypos by\height
\putmorphism(\xpos,\ypos)(-1,-1)[``{#7}]{\height}{\arrowtypea}l%
\putvmorphism(\xpos,\ypos)[#4``{#8}]{\height}{\arrowtypeb}r%
}}

\def\putdtriangle{\@ifnextchar <{\putdtrianglep}{\putdtrianglep
   <\arrowtypea`\arrowtypeb`\arrowtypec;\height>}}
\def\dtriangle{\@ifnextchar <{\dtrianglep}{\dtrianglep
   <\arrowtypea`\arrowtypeb`\arrowtypec;\height>}}
\def\dtrianglep<#1>[#2`#3`#4;#5`#6`#7]{{
\settriparms[#1]
\width=\height                                
\diagram
\putdtrianglep<\arrowtypea`\arrowtypeb`
\arrowtypec;\height>
(0,0)[#2`#3`#4;#5`#6`{#7}]
\enddiagram
}}

\def\putbtrianglep<#1>(#2,#3)[#4`#5`#6;#7`#8`#9]{{%
\settriparms[#1]%
\xpos=#2 \ypos=#3
\puthmorphism(\xpos,\ypos)[#5`#6`{#9}]{\height}{\arrowtypec}b%
\advance\ypos by\height
\putmorphism(\xpos,\ypos)(1,-1)[``{#8}]{\height}{\arrowtypeb}r%
\putvmorphism(\xpos,\ypos)[#4``{#7}]{\height}{\arrowtypea}l%
}}

\def\putbtriangle{\@ifnextchar <{\putbtrianglep}{\putbtrianglep
   <\arrowtypea`\arrowtypeb`\arrowtypec;\height>}}
\def\btriangle{\@ifnextchar <{\btrianglep}{\btrianglep
   <\arrowtypea`\arrowtypeb`\arrowtypec;\height>}}
\def\btrianglep<#1>[#2`#3`#4;#5`#6`#7]{{
\settriparms[#1]
\width=\height                               
\diagram
\putbtrianglep<\arrowtypea`\arrowtypeb`
\arrowtypec;\height>
(0,0)[#2`#3`#4;#5`#6`{#7}]
\enddiagram
}}

\def\putAtrianglep<#1>(#2,#3)[#4`#5`#6;#7`#8`#9]{{%
\settriparms[#1]%
\xpos=#2 \ypos=#3
{\multiply \height by2
\puthmorphism(\xpos,\ypos)[#5`#6`{#9}]{\height}{\arrowtypec}b}%
\advance\xpos by\height \advance\ypos by\height
\putmorphism(\xpos,\ypos)(-1,-1)[#4``{#7}]{\height}{\arrowtypea}l%
\putmorphism(\xpos,\ypos)(1,-1)[``{#8}]{\height}{\arrowtypeb}r%
}}

\def\putAtriangle{\@ifnextchar <{\putAtrianglep}{\putAtrianglep
   <\arrowtypea`\arrowtypeb`\arrowtypec;\height>}}
\def\Atriangle{\@ifnextchar <{\Atrianglep}{\Atrianglep
   <\arrowtypea`\arrowtypeb`\arrowtypec;\height>}}
\def\Atrianglep<#1>[#2`#3`#4;#5`#6`#7]{{
\settriparms[#1]
\width=\height                                     
\diagram
\putAtrianglep<\arrowtypea`\arrowtypeb`
\arrowtypec;\height>
(0,0)[#2`#3`#4;#5`#6`{#7}]
\enddiagram
}}

\def\putAtrianglepairp<#1>(#2)[#3;#4`#5`#6`#7`#8]{{%
\settripairparms[#1]%
\setpos(#2)%
\settokens[#3]%
\puthmorphism(\xpos,\ypos)[\tokenb`\tokenc`{#7}]{\height}{\arrowtyped}b%
\advance\xpos by\height
\puthmorphism(\xpos,\ypos)[\phantom{\tokenc}`\tokend`{#8}]%
{\height}{\arrowtypee}b%
\advance\ypos by\height
\putmorphism(\xpos,\ypos)(-1,-1)[\tokena``{#4}]{\height}{\arrowtypea}l%
\putvmorphism(\xpos,\ypos)[``{#5}]{\height}{\arrowtypeb}m%
\putmorphism(\xpos,\ypos)(1,-1)[``{#6}]{\height}{\arrowtypec}r%
}}

\def\putAtrianglepair{\@ifnextchar <{\putAtrianglepairp}{\putAtrianglepairp%
   <\arrowtypea`\arrowtypeb`\arrowtypec`\arrowtyped`\arrowtypee;\height>}}
\def\Atrianglepair{\@ifnextchar <{\Atrianglepairp}{\Atrianglepairp%
   <\arrowtypea`\arrowtypeb`\arrowtypec`\arrowtyped`\arrowtypee;\height>}}

\def\Atrianglepairp<#1>[#2;#3`#4`#5`#6`#7]{{
\settripairparms[#1]
\settokens[#2]
\width=\height                                
\diagram
\putAtrianglepairp                            
<\arrowtypea`\arrowtypeb`\arrowtypec`
\arrowtyped`\arrowtypee;\height>
(0,0)[{#2};#3`#4`#5`#6`{#7}]
\enddiagram
}}

\def\putVtrianglep<#1>(#2,#3)[#4`#5`#6;#7`#8`#9]{{%
\settriparms[#1]%
\xpos=#2 \ypos=#3
\advance\ypos by\height
{\multiply\height by2
\puthmorphism(\xpos,\ypos)[#4`#5`{#7}]{\height}{\arrowtypea}a}%
\putmorphism(\xpos,\ypos)(1,-1)[`#6`{#8}]{\height}{\arrowtypeb}l%
\advance\xpos by\height
\advance\xpos by\height
\putmorphism(\xpos,\ypos)(-1,-1)[``{#9}]{\height}{\arrowtypec}r%
}}

\def\putVtriangle{\@ifnextchar <{\putVtrianglep}{\putVtrianglep
   <\arrowtypea`\arrowtypeb`\arrowtypec;\height>}}
\def\Vtriangle{\@ifnextchar <{\Vtrianglep}{\Vtrianglep
   <\arrowtypea`\arrowtypeb`\arrowtypec;\height>}}
\def\Vtrianglep<#1>[#2`#3`#4;#5`#6`#7]{{
\settriparms[#1]
\width=\height                                 
\diagram
\putVtrianglep<\arrowtypea`\arrowtypeb`
\arrowtypec;\height>
(0,0)[#2`#3`#4;#5`#6`{#7}]
\enddiagram
}}

\def\putVtrianglepairp<#1>(#2)[#3;#4`#5`#6`#7`#8]{{
\settripairparms[#1]%
\setpos(#2)%
\settokens[#3]%
\advance\ypos by\height
\putmorphism(\xpos,\ypos)(1,-1)[`\tokend`{#6}]{\height}{\arrowtypec}l%
\puthmorphism(\xpos,\ypos)[\tokena`\tokenb`{#4}]{\height}{\arrowtypea}a%
\advance\xpos by\height
\puthmorphism(\xpos,\ypos)[\phantom{\tokenb}`\tokenc`{#5}]%
{\height}{\arrowtypeb}a%
\putvmorphism(\xpos,\ypos)[``{#7}]{\height}{\arrowtyped}m%
\advance\xpos by\height
\putmorphism(\xpos,\ypos)(-1,-1)[``{#8}]{\height}{\arrowtypee}r%
}}

\def\putVtrianglepair{\@ifnextchar <{\putVtrianglepairp}{\putVtrianglepairp%
    <\arrowtypea`\arrowtypeb`\arrowtypec`\arrowtyped`\arrowtypee;\height>}}
\def\Vtrianglepair{\@ifnextchar <{\Vtrianglepairp}{\Vtrianglepairp%
    <\arrowtypea`\arrowtypeb`\arrowtypec`\arrowtyped`\arrowtypee;\height>}}
\def\Vtrianglepairp<#1>[#2;#3`#4`#5`#6`#7]{{
\settripairparms[#1]
\settokens[#2]
\diagram
\putVtrianglepairp                             
<\arrowtypea`\arrowtypeb`\arrowtypec`
\arrowtyped`\arrowtypee;\height>
(0,0)[{#2};#3`#4`#5`#6`{#7}]
\enddiagram
}}

\def\putCtrianglep<#1>(#2,#3)[#4`#5`#6;#7`#8`#9]{{%
\settriparms[#1]%
\xpos=#2 \ypos=#3
\advance\ypos by\height
\putmorphism(\xpos,\ypos)(1,-1)[``{#9}]{\height}{\arrowtypec}l%
\advance\xpos by\height
\advance\ypos by\height
\putmorphism(\xpos,\ypos)(-1,-1)[#4`#5`{#7}]{\height}{\arrowtypea}l%
{\multiply\height by 2
\putvmorphism(\xpos,\ypos)[`#6`{#8}]{\height}{\arrowtypeb}r}%
}}

\def\putCtriangle{\@ifnextchar <{\putCtrianglep}{\putCtrianglep
    <\arrowtypea`\arrowtypeb`\arrowtypec;\height>}}
\def\Ctriangle{\@ifnextchar <{\Ctrianglep}{\Ctrianglep
    <\arrowtypea`\arrowtypeb`\arrowtypec;\height>}}
\def\Ctrianglep<#1>[#2`#3`#4;#5`#6`#7]{{
\settriparms[#1]
\width=\height                               
\diagram
\putCtrianglep<\arrowtypea`\arrowtypeb`
\arrowtypec;\height>
(0,0)[#2`#3`#4;#5`#6`{#7}]
\enddiagram
}}                                           
\def\putDtrianglep<#1>(#2,#3)[#4`#5`#6;#7`#8`#9]{{%
\settriparms[#1]%
\xpos=#2 \ypos=#3
\advance\xpos by\height \advance\ypos by\height
\putmorphism(\xpos,\ypos)(-1,-1)[``{#9}]{\height}{\arrowtypec}r%
\advance\xpos by-\height \advance\ypos by\height
\putmorphism(\xpos,\ypos)(1,-1)[`#5`{#8}]{\height}{\arrowtypeb}r%
{\multiply\height by 2
\putvmorphism(\xpos,\ypos)[#4`#6`{#7}]{\height}{\arrowtypea}l}%
}}

\def\putDtriangle{\@ifnextchar <{\putDtrianglep}{\putDtrianglep
    <\arrowtypea`\arrowtypeb`\arrowtypec;\height>}}
\def\Dtriangle{\@ifnextchar <{\Dtrianglep}{\Dtrianglep
   <\arrowtypea`\arrowtypeb`\arrowtypec;\height>}}
\def\Dtrianglep<#1>[#2`#3`#4;#5`#6`#7]{{
\settriparms[#1]
\width=\height                              
\diagram
\putDtrianglep<\arrowtypea`\arrowtypeb`
\arrowtypec;\height>
(0,0)[#2`#3`#4;#5`#6`{#7}]
\enddiagram
}}                                          
\def\setrecparms[#1`#2]{\width=#1 \height=#2}%

\def\recursep<#1`#2>[#3;#4`#5`#6`#7`#8]{{%
\width=#1 \height=#2
\settokens[#3]
\settowidth{\tempdimen}{$\tokena$}
\ifdim\tempdimen=0pt
  \savebox{\tempboxa}{\hbox{$\tokenb$}}%
  \savebox{\tempboxb}{\hbox{$\tokend$}}%
  \savebox{\tempboxc}{\hbox{$#6$}}%
\else
  \savebox{\tempboxa}{\hbox{$\hbox{$\tokena$}\times\hbox{$\tokenb$}$}}%
  \savebox{\tempboxb}{\hbox{$\hbox{$\tokena$}\times\hbox{$\tokend$}$}}%
  \savebox{\tempboxc}{\hbox{$\hbox{$\tokena$}\times\hbox{$#6$}$}}%
\fi
\ypos=\height
\divide\ypos by 2
\xpos=\ypos
\advance\xpos by \width
\bfig
\putCtrianglep<-1`1`1;\ypos>(0,0)[`\tokenc`;#5`#6`{#7}]%
\puthmorphism(\ypos,0)[\tokend`\usebox{\tempboxb}`{#8}]{\width}{-1}b%
\puthmorphism(\ypos,\height)[\tokenb`\usebox{\tempboxa}`{#4}]{\width}{-1}a%
\advance\ypos by \width
\putvmorphism(\ypos,\height)[``\usebox{\tempboxc}]{\height}1r%
\efig
}}

\def\recurse{\@ifnextchar <{\recursep}{\recursep<\width`\height>}}

\def\puttwohmorphisms(#1,#2)[#3`#4;#5`#6]#7#8#9{{%
\puthmorphism(#1,#2)[#3`#4`]{#7}0a
\ypos=#2
\advance\ypos by 20
\puthmorphism(#1,\ypos)[\phantom{#3}`\phantom{#4}`{#5}]{#7}{#8}a
\advance\ypos by -40
\puthmorphism(#1,\ypos)[\phantom{#3}`\phantom{#4}`{#6}]{#7}{#9}b
}}

\def\puttwovmorphisms(#1,#2)[#3`#4;#5`#6]#7#8#9{{%
\putvmorphism(#1,#2)[#3`#4`]{#7}0a
\xpos=#1
\advance\xpos by -20
\putvmorphism(\xpos,#2)[\phantom{#3}`\phantom{#4}`{#5}]{#7}{#8}l
\advance\xpos by 40
\putvmorphism(\xpos,#2)[\phantom{#3}`\phantom{#4}`{#6}]{#7}{#9}r
}}

\def\puthcoequalizer(#1)[#2`#3`#4;#5`#6`#7]#8#9{{%
\setpos(#1)%
\puttwohmorphisms(\xpos,\ypos)[#2`#3;{#5}`{#6}]{#8}11%
\advance\xpos by #8
\puthmorphism(\xpos,\ypos)[\phantom{#3}`#4`{#7}]{#8}1{#9}
}}

\def\putvcoequalizer(#1)[#2`#3`#4;#5`#6`#7]#8#9{{%
\setpos(#1)%
\puttwovmorphisms(\xpos,\ypos)[#2`#3;{#5}`{#6}]{#8}11%
\advance\ypos by -#8
\putvmorphism(\xpos,\ypos)[\phantom{#3}`#4`#7]{#8}1{#9}
}}

\def\putthreehmorphisms(#1)[#2`#3;#4`#5`#6]#7(#8)#9{{%
\setpos(#1) \settypes(#8)
\if a#9 %
     \vertsize{\tempcounta}{#5}%
     \vertsize{\tempcountb}{#6}%
     \ifnum \tempcounta<\tempcountb \tempcounta=\tempcountb \fi
\else
     \vertsize{\tempcounta}{#4}%
     \vertsize{\tempcountb}{#5}%
     \ifnum \tempcounta<\tempcountb \tempcounta=\tempcountb \fi
\fi
\advance \tempcounta by 60
\puthmorphism(\xpos,\ypos)[#2`#3`{#5}]{#7}{\arrowtypeb}{#9}
\advance\ypos by \tempcounta
\puthmorphism(\xpos,\ypos)[\phantom{#2}`\phantom{#3}`{#4}]{#7}{\arrowtypea}{#9}
\advance\ypos by -\tempcounta \advance\ypos by -\tempcounta
\puthmorphism(\xpos,\ypos)[\phantom{#2}`\phantom{#3}`{#6}]{#7}{\arrowtypec}{#9}
}}

\def\setarrowtoks[#1`#2`#3`#4`#5`#6]{%
\def\toka{#1}
\def\tokb{#2}
\def\tokc{#3}
\def\tokd{#4}
\def\toke{#5}
\def\tokf{#6}
}
\def\hex{\@ifnextchar <{\hexp}{\hexp<1000`400>}}
\def\hexp<#1`#2>[#3`#4`#5`#6`#7`#8;#9]{%
\setarrowtoks[#9]
\yext=#2 \advance \yext by #2
\xext=#1 \advance\xext by \yext
\bfig
\putCtriangle<-1`0`1;#2>(0,0)[`#5`;\tokb``\tokd]
\xext=#1 \yext=#2 \advance \yext by #2
\putsquare<1`0`0`1;\xext`\yext>(#2,0)[#3`#4`#7`#8;\toka```\tokf]
\advance \xext by #2
\putDtriangle<0`1`-1;#2>(\xext,0)[`#6`;`\tokc`\toke]
\efig
}

 \numberwithin{equation}{section}
 \title{Abstract blowing down}
 \author{Michel Van den Bergh}
 \address{Limburgs Universitair Centrum\\ Departement WNI\\ Universitaire
 Campus\\ 3590 Diepenbeek \\ Belgium}
 \email{vdbergh@luc.ac.be, http://www.luc.ac.be/Research/Algebra}
 \thanks{The author is a director of research at the NFWO}
 
\begin{document}
 \begin{abstract}
Assume that $X$ is a  surface over an algebraically closed field
$k$. Let $\tilde{X}$ be obtained from $X$ by blowing up a smooth point and
let $L$ be the exceptional curve. Let
$\coh(X)$ be the category of coherent sheaves on $X$. In this note we
show how to recover $\coh({X})$ from $\coh(\tilde{X})$, if we  know
the object
$\Oscr_L(L)$.  
\end{abstract}
\maketitle
\section{Introduction}
Assume that $X$ is a  surface over an algebraically closed field
$k$. Let $\tilde{X}$ be obtained from $X$ by blowing up a smooth point and
let $L$ be the exceptional curve. Let
$\coh(X)$ be the category of coherent sheaves on $X$. It is an
interesting question to recover 
$\coh(X)$ from $\coh(\tilde{X})$ provided we know $N=\Oscr_L(L)$. 

The reason why we prefer to take $N$ as our basic object
instead of $\Oscr_L$ is the following.
Let $\alpha:\tilde{X}\r X$ be the projection map. Then clearly
$R\alpha_\ast N=0$. So in some sense $\coh(X)$ should be
obtained from $\coh(\tilde{X})$ by dividing out  $N$. 

This does not quite work on the level of abelian categories, but it
does work on the level of derived categories. In \cite{Bondal1} it is
shown that in an appropriate sense
$D^b(\coh(X))=D^b(\coh(\tilde{X}))/[N]$. Hence to find $\coh(X)$ we
have to put a $t$-structure \cite{BBD} on $D^b(\coh(\tilde{X}))/[N]$.
It is easier however to proceed sligthly
differently. From the fact that
$R\alpha_\ast\Oscr_{\tilde{X}}=\Oscr_X$ one deduces that $R\alpha_\ast
L\alpha^\ast$ is the identity on $D^b(\coh(X))$. Hence $L\alpha^\ast $
defines a full and faithful embedding of $D^b(\coh(X))$ in
$D^b(\coh(\tilde{X}))$. Let $\Sscr$ be its essential image. Then
according to \cite{Bondal1} we have
\[
\Sscr=\{A\in D^b(\coh(\tilde{X}))\mid \RHom(A,N)=0
\}
\] 
So we have to put a $t$-structure on $\Sscr$ whose heart is
$\coh(X)$. The solution to this problem is a follows. Define the
following categories 
\begin{align*}
\Tscr&=\{T\in\coh(\tilde{X})\mid \Hom(T,N)=0
\}\\
\Fscr&=\{F\in\coh(\tilde{X})\mid \forall T\in \Tscr: \Hom(T,F)=0
\}
\end{align*}
Following \cite{HRS} we define a ``perverse'' $t$-structure on
$D^b(\coh(\tilde{X}))$ by 
\begin{align*}
\Dscr_{\le 0}^p&=\{B\in\Dscr_{\le 0}\mid H^0(B)\in\Tscr\}\\
\Dscr^p_{\ge 0}&=\{B\in\Dscr_{\ge -1}\mid H^{-1}(B)\in\Fscr\}
\end{align*}
 Then we have our main result.
\begin{theorem} The perverse $t$ structure on $D(\coh(\tilde{X})$
  induces a $t$-structure on $\Sscr$ and
  the heart    of this $t$-stucture is precisely $\coh(X)$.
\end{theorem}
It doesn't seem obvious to  give a version of this result which does not refer to derived categories. 

Below we will give an abstract version of the construction outlined
above. The reason is that exactly the same construction may be
performed in the non-commutative case in order to find an inverse to the
non-commutative blowing up introduced in \cite{VdB19,VdB20}. The details
of this more general case will be published elsewhere.

After this paper was finished I have been informed by  Bondal
that he has independently obtained a similar result (unpublished). I
also wish to thank Bondal for  reading a first version of this
manuscript and for pointing out an error.
\section{Tilting}
In this section we outline a construction from \cite{HRS}.
Let $\Ascr$ be an abelian
category and let $(\Tscr,\Fscr)$ be a \emph{torsion theory} in $\Ascr$.
This is by definition a pair of full additive subcategories of $\Ascr$
such that for every $T\in\Tscr$, $F\in \Fscr$ one has $\Hom(T,F)=0$
and furthermore for every $A$ there exists a (necessarily unique, up
to ismorphism)
exact sequence
\begin{equation}
\label{cond1}
0\r T\r A\r F\r 0
\end{equation}
with $T\in \Tscr$, $F\in \Fscr$. 

From these  conditions it easily follows that 
\begin{align}
\Tscr&=\{T\in\Ascr\mid \forall F\in\Fscr:\Hom(T,F)=0\}\label{cond2}\\
\Fscr&=\{F\in\Ascr\mid \forall T\in\Tscr:\Hom(T,F)=0\}\label{cond3}
\end{align}
Thus $\Tscr={}^\perp \Fscr$ and $\Fscr=\Tscr^\perp$, with obvious
notations. Under suitable finiteness conditions (like $\Ascr$
noetherian) \eqref{cond2}\eqref{cond3} imply \eqref{cond1}

\emph{Tilting} allows one to construct a
new abelian category $\Bscr$ with the roles of $\Tscr$ and $\Fscr$
interchanged. 
Below let $\Dscr$ stand for $D^b(\Ascr)$. One defines
\begin{align*}
\Dscr_{\le 0}^p&=\{B\in\Dscr_{\le 0}\mid H^0(B)\in\Tscr\}\\
\Dscr^p_{\ge 0}&=\{B\in\Dscr_{\ge -1}\mid H^{-1}(B)\in\Fscr\}
\end{align*}
(``$p$'' stands for \emph{perverse}).  It is trivial to verify that
this defines a $t$-structure on $\Dscr$ (using the axioms in \cite[\S
1.3]{BBD}).  Below we will denote the corresponding truncation
functors by $\tau_{\le 0}^p$ and $\tau_{\ge 0}^p$.  

Let
$\Bscr=\Dscr_{\le 0}^p\cap\Dscr^p_{\ge 0}$ be the heart of this
$t$-structure. It follows that this is an abelian category. Let
$(\Tscr',\Fscr')$ be the essential images of $\Fscr[1]$ and $\Tscr$ is
$\Bscr$. Then it is easy to see that $(\Tscr',\Fscr')$ is a torsion
pair in $\Bscr$.
We will call $\Bscr$ the tilting of $\Ascr$ with respect to the
torsion theory $(\Tscr,\Fscr)$. 
\begin{remark} Contrary to what one would expect, it is in
  general not possible to recover $\Ascr$ from $\Bscr$. In order
  for the role of $\Ascr$ and $\Bscr$ to be completely symmetric one
  needs additional conditions. See \cite{HRS}.
\end{remark}
\begin{remark}  If $\Ascr$ is noetherian then this will not in general
  be the case for $\Bscr$. Let $\Ascr$ be the category of finitely
  generated abelian groups and let $\Tscr$, $\Fscr$ be respectively
  the ordinary torsion modules and torsion free modules.
  
  Since $\Ascr$ is hereditary all objects are of the form $F\oplus T$
  for $T\in\Tscr$ and $F\in\Fscr$. Similarly any object in
  $D^b(\Ascr)$ will be the sum of its homology. So every object in
  $\Bscr$ will be of the form $F[1]\oplus T$. It is now a simple
  matter to verify that $\RHom(-,\ZZ[1])$ defines an equivalence
  between $\Bscr$ and $\Ascr^{\text{opp}}$. Since $\Ascr$ is
  noetherian but not artinian, $\Bscr$ will be artinian, but not
  noetherian. 
\end{remark}

\section{The formalism of semi-orthogonal decompositions}
The material in this section is taken from \cite{Bondal2}. 

Let $\Dscr$ be a triangulated category and let $\Bscr$, $\Cscr$ be two
strict full triangulated subcategories of $\Dscr$.  $(\Bscr,\Cscr)$ is said to be a
\emph{semi-orthogonal pair} if $\Hom_\Dscr(B,C)=0$ for $B\in\Bscr$ and $C\in
\Cscr$. Define
\[
\Bscr^\perp=\{A\in\Dscr\mid \forall B\in\Bscr :\Hom_\Dscr(B,A)=0\}
\]
${}^\perp\Cscr$ is defined similarly.

The following result is a slight variation of the statement of
\cite[Lemma 3.1]{Bondal2}.
\begin{lemma} 
\label{bondallemma} The following statements are equivalent for a
  semi-orthogonal pair $(\Bscr,\Cscr)$.
\begin{enumerate}
\item
$\Bscr$ and $\Cscr$ generate $\Dscr$ (as triangulated category).
\item For every $A\in\Dscr$ there exists a distinguished triangle $B\r
  A\r C$ with $B\in\Bscr$ and $C\in\Cscr$.
\item $\Cscr=\Bscr^\perp$ and  the inclusion functor
  $i_\ast:\Bscr\r\Dscr$ has a right adjoint $i^!:\Dscr\r \Bscr$.
\item $\Bscr={}^\perp\Cscr$ and  the inclusion functor $j_\ast:\Cscr\r
  \Dscr$ has a left adjoint  $j^\ast:\Dscr\r \Cscr$
\end{enumerate}
If one of these conditions holds then the triangles in 2. are unique
up to unique isomorphism. They are necessarily of the form 
\begin{equation}
\label{unique}
i_\ast i^! A\r A \r j_\ast j^\ast A
\end{equation}
where the maps are obtained by adjointness from the identiy maps
$i^!A\r i^! A$ and $j^\ast A \r j^\ast A$. In particular triangles as
in 2. are functorial.
\end{lemma}
\begin{remark} The notations $(i_*,i^!, j_*, j^*)$ are purely symbolic
  and shouldn't be interpreted as direct and inverse images. In fact
  in the main application  below  $i_*$ will be
  given by an inverse image!
\end{remark}
If a pair $(\Bscr,\Cscr)$ satisfies one of the conditions of the
previous lemma then we say that it is a \emph{semi-orthogonal
decomposition} of $\Dscr$.
For further reference we note the following diagram of arrows
\begin{equation}
\label{NMfundarrows}
\Cscr \begin{array}{c} j_\ast\\ \rightarrow\\ \leftarrow
  \\ j^\ast\end{array}
 \Dscr 
\begin{array}{c}\displaystyle  i^!\\ \rightarrow\\ 
  \leftarrow \\i_\ast
\end{array}
 \Bscr
\end{equation}
In the following lemma we give some relations between these arrows.
\begin{lemma} One has~:
\begin{gather*}
 i^! i_\ast=\Id_\Bscr
\\
j^\ast j_\ast=\Id_\Cscr\\
j^\ast i_\ast=0\\
i^! j_\ast=0\\
\end{gather*}
\end{lemma}  
In the sequel we will slightly extend the meaning of the notion of
semi-orthogonality. Assume that we have functors
\[
\Cscr\xrightarrow{j_\ast} \Dscr \xleftarrow{i_\ast} \Bscr
\]
which are fully faithful. Assume  furthermore that the essential images of
$\Bscr$ and $\Cscr$ in $\Dscr$ are  semi-orthogonal in
$\Dscr$. Then, if no confusion can arise, we wil also call $(\Bscr,\Cscr)$
is a semi-orthogonal pair in $\Dscr$. Similarly for a semi-orthogonal
decomposition.

\section{Abstract blowing down}
Let $k$ be a field.
and let $\Ascr$ is a noetherian $k$-linear 
abelian category. As before $\Dscr$ stands for $D^b(\Ascr)$. Below we
sometimes need $\RHom(-,-)$ between objects in $D^b(\Ascr)$. In order
to compute this we let $\bar{\Ascr}$ be the closure of $\Ascr$ under
direct limits \cite{Gabriel} and we identify  $D^b(\Ascr)$ with
$D^b_{\Ascr}(\bar{\Ascr})$. We then compute $\RHom$ in the latter
category, which is possible since $\bar{\Ascr}$ has enough
injectives.

Let $N$ be an object of $\Ascr$ satisfying the following
properties. 
\begin{enumerate}
\item
$\RHom(N,N)=k$.
\item For all $A\in\Ascr$ and for all $i$ one has $\dim\Ext^i(A,N)<\infty
$.
\item $\cd \Hom(-,N)\le 2$.
\item The functor $\Ext^2(-,N)^\ast$ is representable.
\end{enumerate}
We will show in this section that it is possible to define a abstract
analogue of the construction outlined in the introduction.

Associated to $N$ we define a torsion theory in $\Ascr$ by
$
\Tscr=\{T\in\Ascr\mid \Hom(T,N)=0\}
$
and $\Fscr=\Tscr^{\perp}$.

The functor $-\Lotimes N$ going from $D^b_f(k)$ to $\Dscr$ has a left
adjoint given by $\RHom(-,N)^\ast$. Furthermore the appropriate
composition of these functors is the identity by condition
1. So $-\Lotimes N$ is a full faithful embedding of $D^b_f(k)$
in $\Dscr$. According to lemma \ref{bondallemma}.4. we can now
construct a diagram as in \eqref{NMfundarrows}
\[
  D^b_f(k)
\bfig
\putmorphism(0,100)(1,0)[``-\Lotimes N]{1100}{1}{a}
\putmorphism(0,-50)(1,0)[``\RHom(-,N)^\ast]{1100}{-1}{b}
\efig
\Dscr
\bfig
\putmorphism(0,100)(1,0)[``R]{300}{1}{a}
\putmorphism(0,-50)(1,0)[``L]{300}{-1}{b}
\efig
\Sscr
\]
where
\[
\Sscr=\{A\in\Dscr\mid \RHom(A,N)=0\}
\]
So $L:\Sscr\r \Dscr$ is the inclusion functor and and $R$ is its
right adjoint.\begin{theorem}
Define
\begin{align*}
\Sscr_{\le 0}&=\Dscr^p_{\le 0} \cap \Sscr\\
\Sscr_{\ge 0}&=\Dscr^p_{\ge 0}\cap \Sscr
\end{align*}
Then $(\Sscr_{\le 0},\Sscr_{\ge 0})$ is a
$t$-structure on $\Sscr$.
\end{theorem}
\begin{proof}
The only non-trivial axiom is  \cite[\S 1.3(iii)]{BBD}
which says that for $X\in\Sscr$ there
should be a triangle $(A,X,B)$ with $A\in\Sscr_{\le 0}$ and $B\in
\Sscr_{\ge 1}$. Now we claim that in fact $\tau^p_{\le 0}X\in \Sscr$ and
$\tau^p_{\ge 1}X\in \Sscr$. This clearly shows what we want.

We need some preparatory work.
\begin{lemma}  If $A\in \Fscr$ then $\Ext^2(A,N)=0$. 
\end{lemma}
\begin{proof} Let $M$ be the object representing $\Ext^2(-,N)^\ast$.
Assume $\Ext^2(A,N)\neq 0$. Then
\[
0\neq \Ext^2(A,N)= \Hom(M,A)^\ast
\]
The proof is complete if we show that $M\in\Tscr$.
However this is clear since
\[
\Hom(M,N)=\Ext^2(N,N)^\ast=0\qed
\]
\def\qed{}\end{proof}
Now apply $\Hom(-,N)$ to the triangle
\[
\tau^p_{\le 0} X \r X \r  \tau^p_{\ge 1} X\r
\]
This yields
\begin{equation}
\label{eeeq1}
\Hom^i(\tau^p_{\le 0} X,N)=\Hom^{i+1}(\tau^p_{\ge 1} X,N)
\end{equation}
It is clear that 
\begin{equation}
\label{eeeq2}
\Hom^i(\tau^p_{\le 0} X,N)=0\qquad \text{for $i\le 0$}
\end{equation}
 and by the previous lemma together with Condition 3. we  also have
\begin{equation}
\label{eeeq3}
\Hom^j(\tau^p_{\ge 1} X,N)=0\qquad \text{for $j\ge 2$}
\end{equation}
Combining \eqref{eeeq1}\eqref{eeeq2}\eqref{eeeq3} yields 
$\Hom^i(\tau^p_{\le 0} X,N)=0$, $\Hom^j(\tau^p_{\ge 1} X,N)=0$ for all
$i,j$. This finishes the proof.
\end{proof}
Let $\Cscr$ be the heart of the $t$-structure on $\Sscr$. We view
$\Cscr$ as the abstract blowing down of $N$ in $\Ascr$. So object in
$\Cscr$ are represented by complexes
\[
A\equiv(A_{-1}\xrightarrow{\theta} A_0)
\]
such that $H_{-1}=\ker \theta\in\Fscr$  and $H_0=\coker
\theta\in\Tscr$  satisfying in addition $\RHom(A,N)=0$.

One problem we have not be able to resolve is the following. 
\begin{question} In the generality above does one necessarily have that $\Cscr$ is noetherian? If not, what additional conditions are necessary?
\end{question}

\section{The commutative case}
What remains to be checked  is that in the commutative our hypotheses are verified, and furthermore that we get the correct answer.

So let us return to the situation from the introduction. So $X$ is a surface over an algebraically closed field $k$ and $\tilde{X}$ is obtained by blowing up a smooth point $p$. So we have a commutative diagram
\[
\begin{CD}
L @>>> \tilde{X}\\
@VVV @V\alpha VV\\
p @>>> X
\end{CD}
\]
where $L$ is the exceptional curve.
We put $\Ascr=\coh(\tilde{X})$, $N=\Oscr_L(L)$. Other notations will be
as above.

To verify the hypotheses we note that they can all be verified in a
neighborhood of $p$. Whence we may assume that $X$ is smooth. As a
next step we will replace $X$ by a smooth compactification. Then
condition 2. is clear and conditions 3,4.  follow from the sophisticated
version of Serre duality
\[
\Ext^i(U,V)=\Ext^{2-i}(V,\omega_{\tilde{X}}\otimes_{\Oscr_{\tilde{X}}} U)^\ast
\]
Condition 1. follows easily from the long exact sequence for $\Hom(-,\Oscr_L(L))$ associated to the short exact sequence
\[
0\r \Oscr_{\tilde{X}}\r \Oscr_{\tilde{X}}(L)\r \Oscr_L(L)\r 0
\]

According to \cite{Bondal1} we  have a semi-orthogonal decomposition
\[
  D^b_f(k)
\bfig
\putmorphism(0,100)(1,0)[``-\Lotimes \Oscr_L(L)]{1100}{1}{a}
\putmorphism(0,-50)(1,0)[``\RHom(-,\Oscr_L(L))^\ast]{1100}{-1}{b}
\efig
D^b(\coh(\tilde{X}))
\bfig
\putmorphism(0,100)(1,0)[``R\alpha_\ast]{300}{1}{a}
\putmorphism(0,-50)(1,0)[``L\alpha^\ast]{300}{-1}{b}
\efig
 D^b(\coh(X))
\]
Thus if we put 
\[
\Sscr=\{A\in \Dscr\mid \RHom(A,N)=0\}
\]
it follows from lemma \ref{bondallemma} that there are inverse equivalences
\[
\Sscr
\bfig
\putmorphism(0,100)(1,0)[``R\alpha_\ast]{300}{1}{a}
\putmorphism(0,-50)(1,0)[``L\alpha^\ast]{300}{-1}{b}
\efig
 D^b_f(\coh(X))
\]
The canonical $t$-structure on $D^b(\coh(X))$ induces a $t$-structure on
$\Sscr$. We have to  show that it coincides with the one we have defined
earlier. That is we have to show that
\begin{align*}
\Sscr_{\le 0}&=L\alpha^\ast \bigl(D^b(\coh(X))_{\le 0}\bigr)\\
\Sscr_{\ge 0}&= L\alpha^\ast \bigl(D^b(\coh(X))_{\ge 0} \bigr)
\end{align*}
Since $\Sscr_{\ge 0}=\Sscr_{\le -1}^\perp$, and similarly for $D^b(\coh(X))_{\ge 0}$ it suffices to verify the first of these equalities. 

Let us first show ``$\subset $". Since $\Sscr_{\le 0}\subset \Sscr$ and
$L\alpha^\ast R\alpha_\ast$ is  the identity on $\Sscr$ it suffices to
show that $R\alpha_\ast\Sscr_{\le 0} \subset D^b(\coh(X))_{\le 0}$. To
this end it suffices to show that if $A\in \Tscr$ then $R^1\alpha_\ast A=0$.

The triangle 
\[
L\alpha^\ast R\alpha_\ast A\r A \r
\RHom(A,N)^\ast\otimes N \r
\]
obtained from \eqref{unique}  yields $H^1(L\alpha^\ast R\alpha_\ast A)=0$. Now $H^1(L\alpha^\ast
R\alpha_\ast A)=\alpha^\ast R^1\alpha_\ast A$. If $R^1\alpha_\ast A\neq
0$ then also $\alpha^\ast R^1\alpha_\ast A\neq 0$ (since $\alpha$ is
surjective) and we obtain a contradiction.

Now we prove the opposite inclusion. Let $T\in D^b(\coh(X))_{\le 0}$.
Clearly $L\alpha^\ast T\in \Sscr\cap \Dscr_{\le 0}$. Furthermore
\[
\Hom(H^0(L\alpha^\ast T),N)=\Hom(L\alpha^\ast T, N)=\Hom(T,R\alpha_\ast
N)=0 \]
whence $H^0(L\alpha^\ast T)\in \Tscr$. This proves what we want.
\ifx\undefined\bysame
\newcommand{\bysame}{\leavevmode\hbox to3em{\hrulefill}\,}
\fi

\end{document}